\RequirePackage{ifpdf}
\ifpdf 
\documentclass[pdftex]{sigma}
\else
\documentclass{sigma}
\fi

\newtheorem{Theorem}{Theorem}[section]
\newtheorem{Lemma}[Theorem]{Lemma}

\newtheorem{Corollary}[Theorem]{Corollary}
{\theoremstyle{definition}

\newtheorem{Definition}[Theorem]{Definition}
\newtheorem{Note}[Theorem]{Note}
}

\begin{document}

\allowdisplaybreaks

\renewcommand{\PaperNumber}{065}

\FirstPageHeading

\ShortArticleName{Double Af\/f\/ine Hecke Algebras of Rank 1}

\ArticleName{Double Af\/f\/ine Hecke Algebras of Rank 1\\ and
the $\boldsymbol{\mathbb Z_3}$-Symmetric Askey--Wilson Relations}

\Author{Tatsuro ITO~$^\dag$ and Paul TERWILLIGER~$^\ddag$}

\AuthorNameForHeading{T. Ito and P. Terwilliger}

\Address{$^\dag$~Division
of Mathematical and Physical Sciences,
Graduate School of Natural Science \\
\hphantom{$^\dag$}~and Technology,
Kanazawa University,
Kakuma-machi,  Kanazawa 920-1192, Japan}
\EmailD{\href{mailto:tatsuro@kenroku.kanazawa-u.ac.jp}{tatsuro@kenroku.kanazawa-u.ac.jp}}

\Address{$^\ddag$~Department
of Mathematics,
University of Wisconsin,\\
\hphantom{$^\ddag$}~480 Lincoln Drive,
Madison, WI 53706-1388, USA}
\EmailD{\href{mailto:terwilli@math.wisc.edu}{terwilli@math.wisc.edu}}

\ArticleDates{Received January 23, 2010, in f\/inal form August 10, 2010;  Published online August 17, 2010}

\Abstract{We consider
the double af\/f\/ine Hecke algebra
$H=H(k_0,k_1,k^\vee_0,k^\vee_1;q)$ associated with
the root system $(C^\vee_1,C_1)$. We display three
elements $x$, $y$, $z$ in $H$ that satisfy
essentially
the $\mathbb Z_3$-symmetric Askey--Wilson relations.
We obtain the relations as follows.
We work with an algebra $\hat H$ that is more general
than $H$,
called the universal double af\/f\/ine Hecke algebra of type
$(C_1^\vee,C_1)$.  An advantage of $\hat H$ over $H$ is
that it is parameter free
and has a larger automorphism group.
We give a surjective algebra homomorphism
${\hat H} \to H$.
We def\/ine some elements $x$, $y$, $z$ in $\hat H$
that get
mapped to their
counterparts in $H$ by this homomorphism.
We give an action of
Artin's braid group $B_3$ on $\hat H$ that acts nicely
on
the elements $x$, $y$, $z$; one generator
sends $x\mapsto y\mapsto z \mapsto x$ and
another generator
interchanges $x$, $y$. Using the $B_3$ action we show that
the elements $x$, $y$, $z$ in $\hat H$ satisfy three equations
that resemble the
$\mathbb Z_3$-symmetric Askey--Wilson relations.
Applying the homomorphism ${\hat H}\to H$
we f\/ind that the elements $x$, $y$, $z$ in $H$ satisfy
similar relations.}


\Classification{33D80; 33D45}

\section{Introduction}\label{section1}

The double af\/f\/ine Hecke algebra (DAHA) for a reduced root
system was def\/ined
by Cherednik~\cite{chered},
 and the def\/inition was extended to include
nonreduced
root systems by
Sahi~\cite{sahi}.
The most ge\-ne\-ral  DAHA of rank 1 is
associated with the root system
 $(C^\vee_1,C_1)$~\cite{oblomkov}; this algebra involves f\/ive nonzero
parameters and will be denoted
by $H=H(k_0,k_1,k^\vee_0,k^\vee_1;q)$.
We mention some recent results
on $H$.
In~\cite{sahi-someprop}
Sahi
links
certain
$H$-modules to
the Askey--Wilson polynomials~\cite{Askey}.
This link is given a comprehensive
treatment
by
Noumi and Stokman~\cite{noumi}.
In~\cite{oblomkov-fd}
Oblomkov
and Stoica
describe the
  f\/inite-dimensional
irreducible $H$-modules under the assumption that
$q$ is not a root of unity.
In~\cite{oblomkov} Oblomkov gives a detailed study
of the algebraic structure of $H$, and f\/inds
an intimate connection to the geometry of
af\/f\/ine cubic surfaces.
His point of departure is the case $q=1$;
under that assumption he f\/inds
that the spherical subalgebra of $H$ is gene\-rated by
three elements $X_1$, $X_2$, $X_3$ that
mutually commute and satisfy a certain cubic equation
\cite[Theo\-rem~2.1, Proposition~3.1]{oblomkov}.
In~\cite{Koo1,Koo2}
 Koornwinder describes the spherical subalgebra of
$H$ under the assumption
 that $q$ is not a root of
unity.
His main results
\cite[Corollary~6.3]{Koo1},
\cite[Theorem~3.2]{Koo2}
 are similar in nature to those of
Oblomkov, although
he formulates these results in a~very
dif\/ferent way and works with
 a dif\/ferent
presentation of $H$.
In Koornwinder's formulation
the spherical subalgebra of $H$ is
related to the Askey--Wilson algebra $AW(3)$, which
was introduced by Zhedanov
in~\cite{Zhidd}.
The original presentation of $AW(3)$ involves
three generators and three relations
\cite[lines (1.1a)--(1.1c)]{Zhidd}.
Koornwinder works with a slightly dif\/ferent
presentation for~$AW(3)$
that involves two generators
and two relations
\cite[lines (2.1), (2.2)]{Koo1}.
These two relations are sometimes called the
Askey--Wilson relations~\cite{aw}.
For the algebra
$AW(3)$
a~third presentation is known
\cite[p.~101]{rosen},
\cite{smithandbell},
\cite[Section 4.3]{wieg}
and
described as follows.
For a sequence of scalars
$g_x$, $g_y$, $g_z$, $h_x$, $h_y$, $h_z$
 the corresponding Askey--Wilson algebra
is def\/ined by generators $X$, $Y$, $Z$ and relations
\begin{gather}
\label{eq:wieg1}
qXY-q^{-1}YX  =  g_zZ+h_z,
\\
\label{eq:wieg2}
qYZ-q^{-1}ZY  =  g_xX+h_x,
\\
\label{eq:wieg3}
qZX-q^{-1}XZ  =  g_yY+h_y.
\end{gather}
We will refer to
(\ref{eq:wieg1})--(\ref{eq:wieg3})
as the
{\it $\mathbb Z_3$-symmetric Askey--Wilson relations}.
Upon eliminating $Z$ in~(\ref{eq:wieg2}),~(\ref{eq:wieg3}) using
(\ref{eq:wieg1}) we obtain
the
Askey--Wilson relations in the variables $X$, $Y$.
Upon substituting $Z'=g_z Z+h_z$ in
(\ref{eq:wieg1})--(\ref{eq:wieg3})
 we recover the the original presentation
 for $AW(3)$ in the variables
 $X$, $Y$, $Z'$.

In this paper we return to the elements
$X_1$, $X_2$, $X_3$ considered by Oblomkov,
although for
 notational convenience we will call them $x$, $y$, $z$.
We show that $x$, $y$, $z$
satisfy
three equations that
 resemble
the $\mathbb Z_3$-symmetric Askey--Wilson relations.
The resemblance is described as follows.
The equations have the
form
(\ref{eq:wieg1})--(\ref{eq:wieg3}) with
$h_x$, $h_y$, $h_z$ not scalars but instead
rational expressions
involving an element $t_1$ that commutes with
each of $x$, $y$, $z$.
The element $t_1$ appears
earlier in the work of Koornwinder
\cite[Def\/inition~6.1]{Koo1}; we will say more
about this at the end of Section~\ref{section2}.
Our derivation of the three equations
is elementary and
 illuminates a role played by Artin's braid
group~$B_3$.

Our proof is summarized  as follows.
Adapting some ideas of Ion and Sahi~\cite{ion}
we work with an algebra $\hat H$ that is more general
than $H$,
called the universal double af\/f\/ine Hecke algebra (UDAHA) of type
$(C_1^\vee,C_1)$.  An advantage of $\hat H$ over $H$ is
that it is parameter free
and has a~larger automorphism group.
We give a surjective algebra homomorphism
${\hat H} \to H$.
We def\/ine some elements $x$, $y$, $z$ in $\hat H$
that get
mapped to their
 counterparts in $H$ by this homomorphism.
Adapting
\cite[Theorem~2.6]{ion}
we give an action of the
braid group $B_3$ on $\hat H$ that acts nicely
on
the elements $x$, $y$, $z$; one generator
sends $x\mapsto y\mapsto z \mapsto x$ and
another generator
interchanges~$x$,~$y$. Using the $B_3$ action we show that
the elements $x$, $y$, $z$ in $\hat H$ satisfy three equations
that resemble
the $\mathbb Z_3$-symmetric Askey--Wilson relations.
Applying the homomorphism ${\hat H}\to H$
we f\/ind that the elements $x$, $y$, $z$ in $H$ satisfy
similar relations.

\section[The double affine Hecke algebra of type $(C^\vee_1,C_1)$]{The double af\/f\/ine Hecke algebra of type $\boldsymbol{(C^\vee_1,C_1)}$}\label{section2}

 Throughout the paper $\mathbb F$ denotes a f\/ield.
An algebra is meant to be associative and have a~1.

We recall the
 double af\/f\/ine
Hecke algebra of type $(C_1^\vee,C_1)$.
For this algebra there are several presentations
in the literature;
one involves three generators
\cite{Koo1,Koo2,sahi2}
and another involves four generators
\cite[p.~160]{mac},
\cite{noumi,
oblomkov,
oblomkov-fd}.
We will use essentially the presentation
of
\cite[p.~160]{mac},
with an adjustment designed to
make explicit the underlying symmetry.

\begin{Definition}
\label{def:dahac}
 Fix  nonzero scalars $ k_0$, $k_1$,
  $k^\vee_0$, $k^\vee_1$, $q$
in  $\mathbb F$.
Let $H(k_0,k_1,k^\vee_0,k^\vee_1;q)$ denote
the $\mathbb F$-algebra  def\/ined by
generators $ t_i$,
 $t^\vee_i$ $(i=0,1)$
 and
relations
\begin{gather}
(t_i-k_i)(t_i-k^{-1}_i) = 0,
\label{eq:daha1}
\\
(t^\vee_i-k^\vee_i)(t^\vee_i-k^{\vee -1}_i) = 0,
\label{eq:daha2}
\\
t^\vee_0 t_0 t^\vee_1 t_1  =  q^{-1}.
\label{eq:daha3}
\end{gather}
This algebra
is called the {\it double affine Hecke algebra}
(or DAHA) of type $(C^\vee_1,C_1)$.
\end{Definition}

\begin{Note}
In
\cite[p.~160]{mac}
Macdonald gives a presentation of
$H$ involving four generators.
To go from his presentation to ours,
multiply each of his generators and the corresponding
 parameter by ${\sqrt{ -1}}$,
and replace his $q$ by $q^2$.
\end{Note}

 The following result is well known; see
for example
\cite[Corollary~1]{sahi2}.

\begin{Lemma}
\label{lem:inv}
Referring to
 Definition~{\rm \ref{def:dahac}}, for $i\in\lbrace 0,1\rbrace$
 the elements
$t_i$, $t^\vee_i$ are invertible and
\begin{gather*}
t_i+t^{-1}_i=k_i+k^{-1}_i,
\qquad
t^\vee_i+t^{\vee -1}_i=k^\vee_i+k^{\vee -1}_i.
\end{gather*}
\end{Lemma}

\begin{proof}
Def\/ine $r_i=k_i+k^{-1}_i-t_i$ and
$r^\vee_i=k^\vee_i+k^{\vee -1}_i-t^\vee_i$.
Using
(\ref{eq:daha1}),
(\ref{eq:daha2})
 we f\/ind
$t_ir_i=r_it_i=1$ and
$t^\vee _ir^\vee_i=r^\vee_i t^\vee_i=1$.
The result follows.
\end{proof}

We now state our main result.
In this result part (ii) follows from
\cite[Theorem 2.1]{oblomkov}; it is included here for
the sake of completeness.

\begin{Theorem}
\label{thm:daha}
In the algebra
$H(k_0,k_1,k^\vee_0,k^\vee_1;q)$
from Definition~{\rm \ref{def:dahac}},
define
\begin{gather*}
x = t^\vee_0 t_1 + \big(t^\vee_0 t_1\big)^{-1},
\qquad
y = t^\vee_1 t_1 + \big(t^\vee_1 t_1\big)^{-1},
\qquad
z = t_0 t_1 + (t_0 t_1)^{-1}.
\label{eq:xyz}
\end{gather*}
Then the following {\rm (i)--(iv)} hold:
\begin{enumerate}\itemsep=0pt
\item[\rm (i)] $t_1$ commutes with each of $x$, $y$, $z$.
\item[\rm (ii)]
Assume
$q^2 =1$. Then $x$, $y$, $z$ mutually commute.
\item[\rm (iii)]
Assume $q^2 \not=1$ and $q^4=1$. Then ${\rm Char}(\mathbb F) \not=2$
and
\begin{gather*}
\frac{xy+yx}{2}  =
\big(k^\vee_0+k^{\vee -1}_0\big)\big(k^\vee_1+k^{\vee -1}_1\big)+
\big(k_0+k^{-1}_0\big)\big(q^{-1}t_1+qt^{-1}_1\big),
\\
\frac{yz+zy}{2}
 =
\big(k^\vee_1+k^{\vee -1}_1\big)\big(k_0+k^{-1}_0\big)+
\big(k^\vee_0+k^{\vee -1}_0\big)\big(q^{-1}t_1+qt^{-1}_1\big),
\\
\frac{zx+xz}{2}  =
\big(k_0+k^{-1}_0\big)\big(k^\vee_0+k^{\vee -1}_0\big)+
\big(k^\vee_1+k^{\vee -1}_1\big)\big(q^{-1}t_1+qt^{-1}_1\big).
\end{gather*}
\item[\rm (iv)]
Assume
$q^4 \not=1$. Then
\begin{gather*}
\frac{qxy-q^{-1}yx}{q^2-q^{-2}}+ z
  =
\frac{(k^\vee_0+k^{\vee -1}_0)(k^\vee_1+k^{\vee -1}_1)+
(k_0+k^{-1}_0)(q^{-1}t_1+qt^{-1}_1)
}{q+q^{-1}}
,
\\
\frac{qyz-q^{-1}zy}{q^2-q^{-2}}+x
 =
\frac{(k^\vee_1+k^{\vee -1}_1)(k_0+k^{-1}_0)+
(k^\vee_0+k^{\vee -1}_0)(q^{-1}t_1+qt^{-1}_1)
}{q+q^{-1}}
,
\\
\frac{qzx-q^{-1}xz}{q^2-q^{-2}}+y
 =
\frac{
(k_0+k^{-1}_0)(k^\vee_0+k^{\vee -1}_0)+
(k^\vee_1+k^{\vee -1}_1)(q^{-1}t_1+qt^{-1}_1)
}{q+q^{-1}}.
\end{gather*}
\end{enumerate}
\end{Theorem}

 The equations in Theorem~\ref{thm:daha}(iv) resemble the
$\mathbb Z_3$-symmetric
Askey--Wilson relations, as we discussed in Section~\ref{section1}.

We will prove Theorem~\ref{thm:daha} in Section~\ref{section5}.

We comment on how Theorem~\ref{thm:daha}
is related to the work of Koornwinder~\cite{Koo1}.
Def\/ine $x$, $y$, $z$ as in Theorem~\ref{thm:daha}.
Then that theorem
describes how
$x$, $y$, $z$, $t_1$ are related.
If we translate
\cite[Def\/inition~6.1, Corollary~6.3]{Koo1}
into the presentation
of Def\/inition~\ref{def:dahac},
then it describes how
$x$,~$y$,~$t_1$ are related, assuming $q$ is not a root of unity
and some constraints on
$k_0$, $k_1$, $k_0^\vee$, $k_1^\vee$.
Under these assumptions and modulo the translation
the following coincide:
(i)
the main relations
\cite[lines (6.2), (6.3)]{Koo1}
of
\cite[Def\/inition~6.1]{Koo1};
(ii) the relations obtained from
the last two equations of
Theorem~\ref{thm:daha}(iv)
by eliminating~$z$ using the f\/irst
equation.

\section[The universal double affine Hecke algebra of type
$(C^\vee_1,C_1)$]{The universal double af\/f\/ine Hecke algebra of type
$\boldsymbol{(C^\vee_1,C_1)}$}\label{section3}

 In our
proof of
Theorem
\ref{thm:daha} we will
initially work with a homomorphic preimage
 $\hat H$ of
$H(k_0,k_1, $ $k^\vee_0,k^\vee_1;q)$
 called the
universal double af\/f\/ine Hecke algebra of type
$(C^\vee_1,C_1)$.
Before we get into the details,
we would like to acknowledge
how  $\hat H$ is related to the work of
Ion and Sahi \cite{ion}.
Given a general DAHA (not just rank 1)
Ion and Sahi construct a group
$\tilde {\mathcal A}$ called the
 double af\/f\/ine Artin group
\cite[Def\/inition~3.4, Theorem~3.10]{ion}.
The given DAHA is a homomorphic image of
the group $\mathbb F$-algebra  $\mathbb F \tilde {\mathcal A}$
\cite[Def\/inition~1.13]{ion}.
For the case
$(C^\vee_1,C_1)$ of the present paper,
their homomorphism has a factorization
${\mathbb F \tilde { \mathcal A}}  \to {\hat H} \to
H(k_0,k_1,k^\vee_0,k^\vee_1;q)$.
In this section and the next
we will obtain some facts about
$\hat H$. We could obtain these facts
from \cite{ion} by applying the
 homomorphism
${\mathbb F \tilde { \mathcal A}}  \to {\hat H}$,
 but for the purpose of
clarity we will prove everything from f\/irst principles.

We now def\/ine
$\hat H$ and describe some of its basic properties.
In Section~\ref{section4} we will discuss how
the group $B_3$ acts on $\hat H$.
In Section~\ref{section5} we
will use the $B_3$ action to prove
Theorem~\ref{thm:daha}.

\begin{Definition}
\label{def:udaha}
Let $\hat H$ denote the $\mathbb F$-algebra
def\/ined by generators $t^{\pm 1}_i$,
 $(t^{\vee}_i)^{\pm 1}$ $(i=0,1)$
and relations
\begin{alignat}{3}
&  t_i t^{-1}_i = t^{-1}_i t_i = 1,
 \qquad &&
 t^\vee_i t^{\vee -1}_i = t^{\vee -1}_i t^\vee_i = 1, &
\label{eq:rel1}
\\
&  t_i+t^{-1}_i \quad {\hbox {\rm is central}},
\qquad &&
 t^\vee_i+t^{\vee -1}_i \quad {\hbox {\rm is central}}, &
\label{eq:rel2}
\\
&  t^\vee_0 t_0 t^\vee_1 t_1 \quad {\hbox {\rm  is central}}. &&&
\label{eq:rel3}
\end{alignat}
We call $\hat H$ the {\it universal double affine Hecke
algebra} (or UDAHA) of type $(C^\vee_1,C_1)$.
\end{Definition}

\begin{Note}
\label{def:artin}
The double af\/f\/ine Artin group
$\tilde {\mathcal A}$ of type
$(C^\vee_1,C_1)$ is def\/ined by
generators
$t^{\pm 1}_i$,
 $(t^{\vee}_i)^{\pm 1}$ $(i=0,1)$
and relations~(\ref{eq:rel1}),~(\ref{eq:rel3})
\cite[Theorem~3.11]{ion}.
\end{Note}

\begin{Definition}
\label{def:qdaha}
Observe that in $\hat H$ the element
$t^\vee_0 t_0 t^\vee_1 t_1$ is invertible; let
$Q$ denote the inverse.
\end{Definition}

\begin{Lemma}
\label{lem:map}
 Given nonzero scalars $ k_0$, $k_1$,
$k^\vee_0$, $k^\vee_1$, $q$
in  $\mathbb F$,
there exists a surjective $\mathbb F$-algebra
homomorphism
${\hat H} \to H(k_0,k_1,k^\vee_0,k^\vee_1;q)$
that sends $Q\mapsto q$ and
$t_i\mapsto t_i, t^\vee_i \mapsto t^\vee_i$ for
$i \in \lbrace 0,1\rbrace $.
\end{Lemma}

\begin{proof}
Compare the def\/ining relations for
$\hat H$ and
$ H(k_0,k_1,k^\vee_0,k^\vee_1;q)$.
\end{proof}

  One advantage of
$\hat H$ over
$H(k_0,k_1,k^\vee_0,k^\vee_1;q)$ is that $\hat H$
has more automorphisms. This is illustrated in the
next lemma.
By an {\it automorphism} of $\hat H$ we mean
an $\mathbb F$-algebra isomorphism ${\hat H} \to {\hat H}$.

\begin{Lemma}
\label{lem:4cycle}
There exists an automorphism
of $\hat H$ that sends
\begin{gather*}
t^\vee_0 \mapsto t_0,
\qquad
t_0 \mapsto t^\vee_1,
\qquad
t^\vee_1 \mapsto t_1,
\qquad
t_1 \mapsto t^\vee_0.
\end{gather*}
This automorphism fixes $Q$.
\end{Lemma}
\begin{proof}
The result follows
from Def\/inition~\ref{def:udaha},
once we verify that
$t_0 t^\vee_1 t_1 t^\vee_0=Q^{-1}$.
This equation holds since each side is equal to
 $t^{\vee -1}_0 Q^{-1} t^\vee_0$.
\end{proof}

\begin{Lemma}
\label{lem:Q}
In the algebra $\hat H$
the element $Q^{-1}$ is equal to each of the following:
\begin{gather}
\label{eq:qversions}
t^\vee_0 t_0 t^\vee_1 t_1,
\qquad
t_0 t^\vee_1 t_1 t^\vee_0,
\qquad
 t^\vee_1 t_1 t^\vee_0 t_0,
\qquad
  t_1 t^\vee_0 t_0 t^\vee_1.
\end{gather}
\end{Lemma}
\begin{proof}
To each side of the equation
$t^\vee_0 t_0 t^\vee_1 t_1 = Q^{-1}$ apply
three times the automorphism
from Lemma~\ref{lem:4cycle}.
\end{proof}

\begin{Definition}
\label{def:xyz}
We def\/ine elements $x$, $y$, $z$ in $\hat H$ as follows.
\begin{gather*}
x = t^\vee_0 t_1 + \big(t^\vee_0 t_1\big)^{-1},
\qquad
y = t^\vee_1 t_1 + \big(t^\vee_1 t_1\big)^{-1},
\qquad
z = t_0 t_1 + (t_0 t_1)^{-1}.
\end{gather*}
\end{Definition}

 The following result suggests why
$x$, $y$, $z$ are of interest.

\begin{Lemma}
\label{lem:why}
Let $u$, $v$ denote invertible elements in any algebra
such that each of $u+u^{-1}$, $v+v^{-1}$ is central.
Then
\begin{enumerate}\itemsep=0pt
\item[\rm (i)] $uv+(uv)^{-1}=vu+(vu)^{-1}$;
\item[\rm (ii)] $uv+(uv)^{-1}$ commutes with each of $u$, $v$.
\end{enumerate}
\end{Lemma}

\begin{proof}
(i) Observe that
\begin{gather*}
uv+(uv)^{-1} =
uv+vu - \big(v+v^{-1}\big)u-v\big(u+u^{-1}\big)+\big(v+v^{-1}\big)\big(u+u^{-1}\big),
\\
vu+(vu)^{-1} =
uv+vu - u\big(v+v^{-1}\big)-\big(u+u^{-1}\big)v+\big(u+u^{-1}\big)\big(v+v^{-1}\big).
\end{gather*}
In these equations
 the expressions on the right
are equal since $u+u^{-1}$ and $v+v^{-1}$ are central.
The result follows.

 (ii) We have
\begin{gather*}
u^{-1}\big(uv+(uv)^{-1}\big)u  =  uv+(uv)^{-1}
\end{gather*}
since each side is equal to $vu+(vu)^{-1}$.
Therefore
$uv+(uv)^{-1}$ commutes with $u$.
 One similarly
shows that
 $uv+(uv)^{-1}$ commutes with $v$.
\end{proof}

\begin{Corollary}
\label{cor:txyz}
In the algebra $\hat H$ the element $t_1$ commutes with each of $x$, $y$, $z$.
\end{Corollary}

\begin{proof}
Use Def\/inition~\ref{def:xyz}
and Lemma~\ref{lem:why}(ii).
\end{proof}

\section[The braid group $B_3$]{The braid group $\boldsymbol{B_3}$}\label{section4}

In this section we display an action of the braid group $B_3$ on
the algebra $\hat H$ from
Def\/inition
\ref{def:udaha}. This $B_3$ action will be used to prove
Theorem
\ref{thm:daha}.

\begin{Definition}
\label{def:braid1}
Artin's braid group $B_3$
is def\/ined by generators $b,c$ and the relation
$b^3=c^2$.
For notational convenience def\/ine
$a=b^3=c^2$.
\end{Definition}

 The following result is a variation
on
\cite[Theorem~2.6]{ion}.

\begin{Lemma}
\label{thm:1} The braid group $B_3$ acts on $\hat H$ as a group
of automorphisms such that $a(h)=t^{-1}_1ht_1$ for
all $h \in {\hat H}$ and $b$, $c$ do the following:
\begin{center}
\begin{tabular}{c|cccc}
$h$  &  $t^\vee_0$ & $t_0$ & $t^\vee_1$ & $t_1$  \\
\hline
$b(h)$ &
$t^{-1}_1t^\vee_1 t_1$ & $t^\vee_0$ & $t_0$ & $t_1$\tsep{2pt}\bsep{2pt}
\\
$c(h)$ &
$t^{-1}_1t^\vee_1 t_1$ & $t^\vee_0 t_0 t^{\vee -1}_0$ & $t^\vee_0$ & $t_1$
\end{tabular}
\end{center}
\end{Lemma}

\begin{proof}
There exists an automorphism $A$ of $\hat H$ that
sends $h\mapsto t^{-1}_1h t_1$ for all $h \in {\hat H}$.
Def\/ine
\begin{gather}
\label{eq:capt}
T^\vee_0=t^{-1}_1 t^{\vee}_1 t_1,
\qquad
T_0 = t^\vee_0,
\qquad
T^\vee_1 = t_0,
\qquad
T_1 = t_1.
\end{gather}
Note that
$T^\vee_0$,
$T_0$,
$T^\vee_1$,
$T_1$
are invertible and that
\begin{alignat*}{3}
& T^\vee_0+T^{\vee -1}_0=t^{\vee}_1 + t^{\vee -1}_1,
\qquad &&
T_0 + T^{-1}_0 = t^\vee_0 + t^{\vee -1}_0,&
\\
& T^\vee_1 + T^{\vee -1}_1= t_0+t^{-1}_0,
\qquad &&
T_1 +T^{-1}_1 = t_1+t^{-1}_1. &
\end{alignat*}
In each of these four equations the expression on the right is central
so the expression on the left is central.
Using~(\ref{eq:capt}) and
Lemma~\ref{lem:Q},
\begin{gather*}
T^\vee_0
T_0
T^\vee_1
T_1 = t^{-1}_1 t^\vee_1 t_1 t^\vee_0 t_0 t_1
 = t^{-1}_1 Q^{-1} t_1 = Q^{-1}
\end{gather*}
so
$T^\vee_0
T_0
T^\vee_1
T_1$ is central.
By these comments there exists an $\mathbb F$-algebra
homomorphism $B:{\hat H}\to {\hat H}$ that
sends
\begin{gather*}
t^\vee_0 \mapsto T^\vee_0,
\qquad
t_0 \mapsto T_0,
\qquad
t^\vee_1 \mapsto T^\vee_1,
\qquad
t_1 \mapsto T_1.
\end{gather*}
We claim that $B^3=A$. To prove the claim we show that
$B^3$, $A$ agree at each of $t^\vee_0$, $t_0$,
$t^\vee_1$, $t_1$.
Note that $A$ f\/ixes $t_1$. Note also
that $t_1$ is f\/ixed by $B$ and hence $B^3$;
therefore
$B^3$ and $A$ agree at~$t_1$.
 The map $B$ sends
\begin{gather*}
t^\vee_1 \mapsto t_0 \mapsto t^\vee_0 \mapsto
 t^{-1}_1t^\vee_1 t_1
\mapsto
 t^{-1}_1t_0 t_1
\mapsto
 t^{-1}_1t^\vee_0 t_1.
\end{gather*}
Therefore $B^3$ sends
\begin{gather*}
t^\vee_1 \mapsto t^{-1}_1 t^\vee_1 t_1,
\qquad
t_0 \mapsto t^{-1}_1 t_0 t_1,
\qquad
t^\vee_0 \mapsto t^{-1}_1 t^\vee_0 t_1,
\end{gather*}
so $B^3$, $A$ agree at each of
$t^\vee_1$, $t_0$, $t^\vee_0$.
We have shown $B^3=A$. By this and since $A$ is invertible,
we see that $B$ is invertible and hence an
automorphism of~$\hat H$.
Def\/ine
\begin{gather}
\label{eq:scap}
S^\vee_0=t^{-1}_1 t^{\vee}_1 t_1,
\qquad
S_0 = t^\vee_0t_0t^{\vee -1}_0,
\qquad
S^\vee_1 = t^\vee_0,
\qquad
S_1 = t_1.
\end{gather}
Note that
$S^\vee_0$,
$S_0$,
$S^\vee_1$,
$S_1$
are invertible and
\begin{alignat*}{3}
& S^\vee_0+S^{\vee -1}_0=t^{\vee}_1 + t^{\vee -1}_1,
\qquad &&
S_0 + S^{-1}_0 = t_0 + t^{-1}_0,&
\\
& S^\vee_1 + S^{\vee -1}_1= t^\vee_0+t^{\vee -1}_0,
\qquad &&
S_1 +S^{-1}_1 = t_1+t^{-1}_1.
\end{alignat*}
In each of these four equations the expression on the right is central
so the expression on the left is central.
Using~(\ref{eq:scap}) and
Lemma~\ref{lem:Q},
\begin{gather*}
S^\vee_0
S_0
S^\vee_1
S_1 = t^{-1}_1 t^\vee_1 t_1 t^\vee_0 t_0
t_1 = t^{-1}_1 Q^{-1} t_1 = Q^{-1}
\end{gather*}
so
$S^\vee_0
S_0
S^\vee_1
S_1$ is central.
By these comments there exists an $\mathbb F$-algebra
homomorphism $C:{\hat H}{\to} {\hat H}$ that
sends
\begin{gather*}
t^\vee_0 \mapsto S^\vee_0,
\qquad
t_0 \mapsto S_0,
\qquad
t^\vee_1 \mapsto S^\vee_1,
\qquad
t_1 \mapsto S_1.
\end{gather*}
We claim that $C^2=A$. To prove the claim we show that
$C^2$, $A$ agree at each of $t^\vee_0$, $t_0$,
$t^\vee_1$, $t_1$. Both
$C^2$ and $A$ f\/ix $t_1$.
 The map $C$ sends
$t^\vee_0 \mapsto t^{-1}_1 t^\vee_1 t_1
\mapsto
t^{-1}_1 t^\vee_0 t_1$ so $C^2$, $A$ agree at
$t^\vee_0$.
 The map $C$ sends
$t^\vee_1 \mapsto t^\vee_0 \mapsto t^{-1}_1 t^\vee_1 t_1$
so $C^2$, $A$ agree at
$t^\vee_1$.
The map $C$ sends
\begin{gather*}
t_0 \mapsto t^\vee_0t_0t^{\vee -1}_0 \mapsto
t^{-1}_1 t^\vee_1t_1t^\vee_0 t_0 t^{\vee -1}_0t^{-1}_1 t^{\vee -1}_1 t_1.
\end{gather*}
In the above line the expression on the right equals $t^{-1}_1t_0t_1$.
To see this, note that $t^\vee_1 t_1 t^\vee_0 t_0 =
t_0 t^\vee_1 t_1 t^\vee_0$ since each side equals~$Q^{-1}$
by Lemma~\ref{lem:Q}.
We have shown that $C^2$, $A$ agree at~$t_0$. By the above comments
$C^2$, $A$ agree at each of
$t^\vee_0$, $t_0$,
$t^\vee_1$, $t_1$ so
$C^2=A$. Therefore $C$ is invertible
and hence an automorphism of $\hat H$.
We have shown that the desired $B_3$ action exists.
\end{proof}

  The next result is
immediate from
 Lemma~\ref{thm:1} and its proof.

\begin{Lemma}
\label{thm:half}
The $B_3$ action from Lemma~{\rm \ref{thm:1}} does the following
to the central elements~\eqref{eq:rel2},~\eqref{eq:rel3}. The generator $a$ fixes every central
element. The generators $b$, $c$ fix $Q$ and satisfy
the table below.
\begin{center}
\begin{tabular}{c|cccc}
$h$  &  $t^\vee_0+t^{\vee -1}_0$ & $t_0+t^{-1}_0$ & $t^\vee_1+t^{\vee -1}_1$
 & $t_1+t^{-1}_1$  \\
\hline
$b(h)$ &
$t^\vee_1 +t^{\vee -1}_1$ & $t^\vee_0+t^{\vee -1}_0$ & $t_0+t^{-1}_0$ & $t_1+t^{-1}_1$\tsep{2pt}\tsep{2pt}
\\
$c(h)$ &
$t^\vee_1+t^{\vee -1}_1$ & $t_0 +t^{-1}_0$ & $t^\vee_0+t^{\vee -1}_0$ & $t_1+t^{-1}_1$
\end{tabular}
\end{center}
\end{Lemma}

\section{The proof of Theorem \ref{thm:daha}}\label{section5}

Recall the elements
$x$, $y$, $z$ of $\hat H$ from
Def\/inition~\ref{def:xyz}.
 In this section
we describe how
the group $B_3$
acts on these elements.
Using this information we show that $x$, $y$, $z$
satisfy three equations that resemble the
$\mathbb Z_3$-symmetric Askey--Wilson relations.
Using these equations we obtain
Theorem~\ref{thm:daha}.

\begin{Theorem}
\label{thm:2}
The $B_3$ action from Lemma~{\rm \ref{thm:1}} does the following
to the elements $x$, $y$, $z$ from
Definition~{\rm \ref{def:xyz}}.
The generator $a$ fixes each of $x$, $y$, $z$. The
generator $b$ sends $x\mapsto y \mapsto z \mapsto x$.
The generator $c$ swaps $x$, $y$ and sends $z \mapsto z'$ where
\begin{gather*}
 Q z + Q^{-1}z' + xy  =
Q^{-1} z + Qz' + yx
\\
 \phantom{Q z + Q^{-1}z' + xy }{}=
 \big(t^\vee_0+t^{\vee -1}_0\big)\big(t^\vee_1+t^{\vee -1}_1\big)+
\big(t_0+t^{-1}_0\big)\big(Q^{-1}t_1+Qt^{-1}_1\big).
\end{gather*}
\end{Theorem}

\begin{proof}
The generator $a$ f\/ixes each of $x$, $y$, $z$ by
Corollary~\ref{cor:txyz} and since $a(h)=t^{-1}_1 h t_1$ for
all $h \in {\hat H}$.
The generator $b$ sends $x\mapsto y \mapsto z \mapsto x$
by Def\/inition~\ref{def:xyz}, Corollary~\ref{cor:txyz},
and Lemma~\ref{thm:1}.
Similarly the generator
$c$ swaps $x$, $y$. Def\/ine $z'=c(z)$.
We show that $z'$ satisf\/ies the equations in
the theorem statement.
We f\/irst show that
\begin{gather}
 Q^{-1} t_0 +
Q c(t_0) + y t^\vee_0 =
(t^\vee_1 t_1)^{-1}\big(t^\vee_0+t^{\vee -1}_0\big) + Q^{-1}\big(t_0+t^{-1}_0\big).
\label{eq:sigy1}
\end{gather}
By Lemma~\ref{thm:1},
$c(t_0) = t^\vee_0 t_0 t^{\vee -1}_0$. By
this and
Def\/inition~\ref{def:qdaha},
\begin{gather}
\label{eq:p1}
Q c(t_0)  =  \big(t^\vee_1 t_1\big)^{-1} t^{\vee -1}_0.
\end{gather}
By Lemma~\ref{lem:Q},
\begin{gather}
\label{eq:p2}
t^\vee_1 t_1 t^\vee_0 = Q^{-1}t^{-1}_0.
\end{gather}
Using~(\ref{eq:p1}),
(\ref{eq:p2}) and
$y=t^\vee_1 t_1+
(t^\vee_1 t_1)^{-1}$ we obtain~(\ref{eq:sigy1}).
Next we show that
\begin{gather}
Q^{-1} t^{-1}_0 +
Q c\big(t^{-1}_0\big) + y t^{\vee -1}_0 =
t^\vee_1 t_1\big(t^\vee_0+t^{\vee -1}_0\big) + Q\big(t_0+t^{-1}_0\big).
\label{eq:sigy2}
\end{gather}
By Lemma~\ref{thm:half},
\begin{gather*}
c(t_0)+
c\big(t^{-1}_0\big) =
 t_0+t^{-1}_0.
\end{gather*}
Combining this with~(\ref{eq:sigy1})
we obtain~(\ref{eq:sigy2}) after a brief calculation.
In~(\ref{eq:sigy1})
we multiply each term on the right by
$t_1$ and use $c(t_1)=t_1$ to get
\begin{gather}
\label{eq:first}
Q^{-1} t_0t_1 +
Qc(t_0 t_1)+ y t^\vee_0 t_1  =
\big(t^\vee_1t_1\big)^{-1} t_1 \big(t^\vee_0 + t^{\vee -1}_0\big)
+
Q^{-1}t_1
 \big(t_0+t^{-1}_0\big).
\end{gather}
In~(\ref{eq:sigy2})
we multiply each term on the left by
$t^{-1}_1$ and use $c(t^{-1}_1)=t^{-1}_1$ together
with the fact that $y$ commutes with $t_1$ to get
\begin{gather}
\label{eq:second}
 Q^{-1} (t_0t_1)^{-1}
+
Qc\big((t_0 t_1)^{-1}\big)+ y \big(t^\vee_0 t_1\big)^{-1}
=
t^{-1}_1 t^\vee_1t_1 \big(t^\vee_0 + t^{\vee -1}_0\big)
+
Qt^{-1}_1
 \big(t_0+t^{-1}_0\big).
\end{gather}
We have
\begin{gather}
\big(t^\vee_1 t_1\big)^{-1} t_1 + t^{-1}_1 t^\vee_1 t_1 = t^\vee_1+
t^{\vee -1}_1
\label{eq:third}
\end{gather}
since both sides equal
$t^{-1}_1(t^\vee_1+t^{\vee -1}_1)t_1$.
We now add (\ref{eq:first}), (\ref{eq:second})
and simplify the result using~(\ref{eq:third}) to obtain
\begin{gather}
Q^{-1}z+ Q z' + yx  =
\big(t^\vee_0+ t^{\vee -1}_0\big)
\big(t^\vee_1+ t^{\vee -1}_1\big)
+
\big(t_0+ t^{-1}_0\big)
\big(Q^{-1}t_1+ Qt^{-1}_1\big).
\label{eq:Qz}
\end{gather}
We now apply~$c$ to each side of
(\ref{eq:Qz}) and evaluate the result.
To aid in this evaluation we
recall that $c$ swaps $x$, $y$; also
 $c$ swaps $z$, $z'$ since
 $c^2=a$ and $a(z)=z$.
By these comments and
Lemma~\ref{thm:half} we obtain
\begin{gather*}
Q z + Q^{-1}z' + xy
=\big(t^\vee_0+t^{\vee -1}_0\big)\big(t^\vee_1+t^{\vee -1}_1\big)+
\big(t_0+t^{-1}_0\big)\big(Q^{-1}t_1+Qt^{-1}_1\big).\tag*{\qed}
\end{gather*}
\renewcommand{\qed}{}
\end{proof}

\begin{Theorem}
\label{thm:main}
In the algebra $\hat H$ the elements $x$, $y$, $z$
are related as follows:
\begin{gather*}
  Qxy-Q^{-1}yx +\big(Q^2-Q^{-2}\big)z
\\
  \qquad  {}= \big(Q-Q^{-1}\big)
\big( \big(t^\vee_0+t^{\vee -1}_0\big)\big(t^\vee_1+t^{\vee -1}_1\big)+
\big(t_0+t^{-1}_0\big)\big(Q^{-1}t_1+Qt^{-1}_1\big) \big),
\\
Qyz-Q^{-1}zy +\big(Q^2-Q^{-2}\big)x
\\
 \qquad {} = \big(Q-Q^{-1}\big)
\big(
\big(t^\vee_1+t^{\vee -1}_1\big)\big(t_0+t^{-1}_0\big)+
\big(t^\vee_0+t^{\vee -1}_0\big)\big(Q^{-1}t_1+Qt^{-1}_1\big) \bigr),
\\
  Qzx-Q^{-1}xz +\big(Q^2-Q^{-2}\big)y
\\
   \qquad
{} = \big(Q-Q^{-1}\big)
\big(
\big(t_0+t^{-1}_0\big)\big(t^\vee_0+t^{\vee -1}_0\big)+
\big(t^\vee_1+t^{\vee -1}_1\big)\big(Q^{-1}t_1+Qt^{-1}_1\big)
\bigr).
\end{gather*}
\end{Theorem}

\begin{proof}
To get the f\/irst equation, eliminate $z'$ from the
equations of Theorem~\ref{thm:2}. To get the other two equations
use the $B_3$ action
from Lemma~\ref{thm:1}.
Specif\/ically,
apply $b$ twice to the f\/irst equation
and use the data in
Lemma~\ref{thm:half}, together with the
fact that $b$
cyclically permutes $x$,~$y$,~$z$.
\end{proof}

\begin{proof}[Proof of Theorem~\ref{thm:daha}]
 Apply the homomorphism ${\hat H} \to H(k_0, k_1, k^\vee_0, k^\vee_1)$
from Lemma
\ref{lem:map}. Part (i) follows via Corollary
\ref{cor:txyz}, and parts (ii)--(iv) follow from
Theorem~\ref{thm:main} together with Lemma~\ref{lem:inv}.
\end{proof}

\subsection*{Acknowledgements}
We thank Alexei Zhedanov for mentioning
to us around 2005 that $AW(3)$ has the presentation
(\ref{eq:wieg1})--(\ref{eq:wieg3}); this knowledge
motivated us to search for a result like
Theorem~\ref{thm:daha}.
 We also thank
Zhedanov for several illuminating conversations on  DAHA
during his visit to Kanazawa in December 2007.
We thank the two referees for clarifying how the
present paper is related to the previous literature.
The second author thanks
Tom Koornwinder,
 Alexei Oblomkov, and
Xiaoguang  Ma
 for useful recent
conversations on the general subject of DAHA.

\pdfbookmark[1]{References}{ref}
\LastPageEnding

\end{document}